\documentclass[english]{article}
\usepackage[T1]{fontenc}
\usepackage[latin9]{inputenc}
\usepackage{amsmath}
\usepackage{graphicx}
\usepackage{amssymb}

\makeatletter

\providecommand{\tabularnewline}{\\}

\makeatother

\usepackage{babel}

\begin{document}

\title{An Algorithm on Multivariate Interpolation for Cartesian Grid Complex
Boundary Problems}

\author{Shuqiang Wang}

\date{06/20/2013}

\maketitle
\tableofcontents{}
\begin{abstract}
We present a simple algorithm to select multivariate interpolation
stencil with a Cartesian grid. We show its applicability by using
this algorithm in the embedded boundary method for solving the elliptic
interface problem.
\end{abstract}

\section{Introduction}

In order to solve Cartesian grid complex boundary problems, we frequently
need to select stencil from the Cartesian grid points for multivariate
interpolation to approximate some function or its derivatives. Generally,
second degree polynomial would be accurate enough for many applications.
In two dimension ($2D$), an arbitrary second degree polynomial has
six unknown coefficients and it is comparatively easy to find a stencil
with six Cartesian grid points which can uniquely determine the polynomial.
However, it is not so easy to do so in three dimension (3D) for a
second degree polynomial with ten unknown coefficients which requires
ten grid points. 

To approximate solution flux (or function derivative), Johansen et
al \cite{JohansenColella98} uses seven grid points to compute a second
order flux in $2D$, and Schwartz et al \cite{Colella2006} uses $19$
Cartesian grid points to compute a second order flux in $3D$. When
we use these fluxes to approximate the elliptic problem, the matrix
has large non-zero entry and thus requires more time to obtain convergence.
Another difficulty is that there might not be enough grid points available
to select from for complex boundary.

Guenther et al \cite{GuentherRoetman1970} gives a simple and clever
way to select stencil for interpolation in higher dimension. Gasca
et al \cite{GascaSauer2000} reviews the polynomial interpolation
algorithms in several variables. With some modification of the idea
in \cite{GuentherRoetman1970}, we present a simple algorithm to select
stencil for multivariate interpolation in $2D$ and $3D$. 

In section \ref{sec:Methods}, we give the algorithms to select the
stencil and also the rationale why the stencil selected can be used
to uniquely determine an arbitrary second degree polynomial. In section
\ref{sec:Examples}, we apply this algorithm in our embedded boundary
method (EBM) \cite{WangSamulyakGuo2010} to solve an elliptic interface
problem. Finally, we give the conclusion.

\section{Methods\label{sec:Methods}}

\subsection{2D Interpolation}

\begin{figure}

\caption{\label{fig:2D-Stencil}2D Stencil}

\begin{centering}
\includegraphics[scale=0.3]{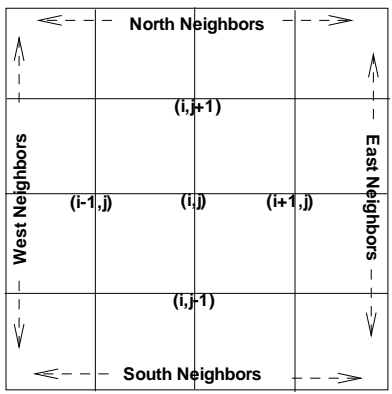}
\par\end{centering}

\end{figure}

Our $2D$ stencil selection algorithm is a modification of the algorithm
presented in \cite{WangSamulyakGuo2010}. For simplicity, we only
select grid points which are at most two grid block away from the
starting grid point $(i,j)$ as shown in Fig \ref{fig:2D-Stencil}.
As in Wang et al\cite{WangSamulyakGuo2010}, we classify the candidate
grid points into different group of neighbors.
\begin{itemize}
\item Zeroth layer neighbor: simply the grid point $(i,j)$ around which
we want to select candidate grid points for interpolation.
\item First layer neighbor: grid points with grid index $(\hat{i},\hat{j})$
satisfying \[
max(\left|\hat{i}-i\right|,\left|\hat{j}-j\right|)=1.\]

\end{itemize}
Thus, the first layer grid points include $(i-1,j-1)$, $(i,j-1)$,
$(i+1,j-1)$, $(i-1,j)$, $(i+1,j)$, $(i-1,j+1)$, $(i,j+1)$. 
\begin{itemize}
\item Second layer neighbor: grid points with grid index $(\hat{i},\hat{j})$
satisfying\[
max(\left|\hat{i}-i\right|,\left|\hat{j}-j\right|)=2.\]

\end{itemize}
For convenience, we separate the second layer neighbors into four
overlapping sets as in Fig \ref{fig:2D-Stencil}: west, east, south,
and north neighbors.
\begin{itemize}
\item Consecutive grid points: consecutive grid points in a given layer.
\end{itemize}
To find 6 points for quadratic polynomial interpolation, Wang et al\cite{WangSamulyakGuo2010}
uses the following algorithm: 1) select grid point $(i,j)$; 2) select
two consecutive grid points from the first layer neighbors; 3) select
three consecutive grid points from the second layer neighbors. However,
it is later realized that one configuration (and its symmetric configurations)
gives a singular matrix when it is used to determine the unknown coefficients
for a quadratic polynomial as in Fig 2.

Our improved algorithm to select candidate cells for quadratic polynomial
interpolation in $2D$ is the following:
\begin{enumerate}
\item always select grid point $(i,j$).
\item select two consecutive grid points from the first layer neighbors.
\item select one direction, say south, and select three grid points from
the second layer south neighbors.
\end{enumerate}
Note that the first two steps of the above algorithm guarantee that
the three selected points are not on a straight line. 

In the following, we prove that we can find a unique quadratic polynomial
using the above algorithm. 

\emph{Proof}. For simplicity, we assume that we select south direction
in the third step. An arbitrary $2D$ quadratic polynomial can be
written as \[
f(x,y)=a_{0}+a_{10}x+a_{01}y+a_{20}x^{2}+a_{11}xy+a_{02}y^{2}\]
which can also be rewritten as \begin{equation}
f(x,y)=a_{0}+a_{10}x+a_{20}x^{2}+y(a_{01}+a_{11}x+a_{02}y).\label{eq:quadratic-polynomial-2d}\end{equation}
We put the origin of the $2D$ coordinate system at the grid point
$(i,j-2)$. Thus the south neighbor grid points are on the $x$ coordinate.
Now we need to determine the unknown coefficients in equation \ref{eq:quadratic-polynomial-2d}.
Note that \[
f(x,0)=a_{0}+a_{10}x+a_{20}x^{2}.\]
There are three unknown coefficients. We can select three grid points
from the second layer south neighbors. After $a_{0}$, $a_{10}$,
and $a_{20}$ are determined, we have \[
f(x,y)-(a_{0}+a_{10}x+a_{20}x^{2})=y(a_{01}+a_{11}x+a_{02}y).\]
When $y\neq0$, we have \[
\left.\frac{f(x,y)-(a_{0}+a_{10}x+a_{20}x^{2})}{y}\right|_{y\neq0}=a_{01}+a_{11}x+a_{02}y.\]
The right hand side of the above equation is a first deegree linear
polynomial with three unknown coefficients. It can be uniquely determined
by using the three grid points selected by the first two steps of
the algorithm since these three points are not on a straight line.
Thus we have proved that the algorithm can uniquely determine a quadratic
polynomial in $2D$. $\blacksquare$

Note that we can easily extend the above proof to modify our algorithm
to find stencil for higher order polynomials. An arbitrary third degree
polynomial can be written as \begin{eqnarray*}
f(x,y) & = & a_{0}+a_{10}x+a_{01}y+a_{20}x^{2}+a_{11}xy+a_{02}y^{2}+\\
 &  & a_{30}x^{3}+a_{21}x^{2}y+a_{12}xy^{2}+a_{03}y^{3}\end{eqnarray*}
which can be rearranged as \begin{eqnarray*}
f(x,y) & = & a_{0}+a_{10}x+a_{20}x^{2}+a_{30}x^{3}+\\
 &  & y(a_{01}+a_{11}x+a_{21}x^{2}+a_{02}y+a_{12}xy+a_{03}y^{2}).\end{eqnarray*}
Suppose that we have select the south direction to select candidate
grid points on the third layer south neighbors. As before, we can
set the origin of the coordinate system at the center of the third
layer south neighbors. Selecting four cells from the third layer south
neighbors, we can uniquely determine $a_{0}$, $a_{10}$, $a_{20}$,
and $a_{30}$. Then we have \[
\left.\frac{f(x,y)-(a_{0}+a_{10}x+a_{20}x^{2}+a_{30}x^{3})}{y}\right|_{y\neq0}=(a_{01}+a_{11}x+a_{21}x^{2}+a_{02}y+a_{12}xy+a_{03}y^{2})\]
where the right hand side is a quadratic polynomial and it can be
determined using the algorithm we have presented above using the grid
point $(i,j)$ and its first and second layer neighbors.

.%
\begin{figure}

\caption{\label{fig:Stencil_2D_singularMatrixConfiguration}2D Stencil with
configuration leading to a singular matrix. The selected grid points
are denoted using a circle at the cell center}

\begin{centering}
\includegraphics[scale=0.3]{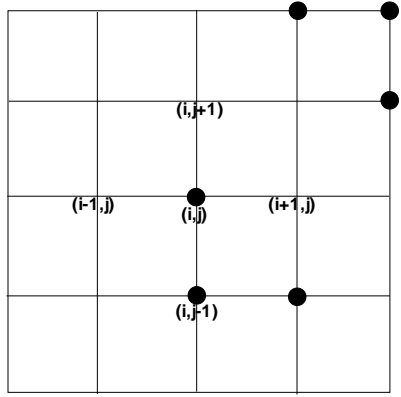}
\par\end{centering}

\end{figure}

\subsection{3D Interpolation}

For a grid point $(i,j,k$), it has 6 directions: west, east in the
x coordinate direction; south, north in the y coordinate direction;
and down, up in the z coordinate direction.

The algorithm for $2D$ can be easily extended to $3D$ as the following:
\begin{enumerate}
\item always select grid point $(i,j,k$).
\item select one direction, say east, and select three consecutive grid
points which are not on a straight line from the first layer east
neighbors.
\item select one direction, say down, and select six grid points from the
second layer down neighbors.
\end{enumerate}
Note that the first two steps of the above algorithm guarantee that
the four selected grid points are not on a plane. Step 2 can be done
using the first two steps of the $2D$ algorithm when we consider
the first layer east neighbors as a $2D$ plane. Similarly, step 3
above can be done using the $2D$ algorithm when we consider the second
layer down neighbors as a $2D$ plane.

Now we prove that the above algorithm can be used to select stencil
to uniquely determine a quadratic polynomial uniquely in $3D$.

\emph{Proof}. The proof is similar to the proof in $2D$. An arbitrary
$3D$ quadratic polynomial can be written as \begin{eqnarray*}
f(x,y,z) & = & a_{0}+a_{100}x+a_{010}y+a_{001}z+\\
 &  & a_{200}x^{2}+a_{110}xy+a_{101}xz+\\
 &  & a_{020}y^{2}+a_{011}yz+a_{002}z^{2}\end{eqnarray*}
which can also be rearranged as \begin{eqnarray*}
f(x,y,z) & = & a_{0}+a_{10}x+a_{01}y+a_{20}x^{2}+a_{11}xy+a_{02}y^{2}+\\
 &  & z(a_{001}+a_{101}x+a_{011}y+a_{002}z)\end{eqnarray*}
Assuming that we select the down direction to select six grid points
from the second layer neighbors. Set the origin at grid point $(i,j,k-2)$.
When $z=0$, we have \begin{eqnarray*}
f(x,y,0) & = & a_{0}+a_{10}x+a_{01}y+a_{20}x^{2}+a_{11}xy+a_{02}y^{2}\end{eqnarray*}
where the unknown coefficients can be determined using the $2D$ algorithm
to select six grid points in the second layer down neighbors. After
$a_{0}$, $a_{10}$, $a_{01}$, $a_{20}$, $a_{11}$, and $a_{02}$
have been determined, we have \[
\left.\frac{f(x,y,z)-a_{0}+a_{10}x+a_{01}y+a_{20}x^{2}+a_{11}xy+a_{02}y^{2}}{z}\right|_{z\neq0}=a_{001}+a_{101}x+a_{011}y+a_{002}z\]
Now, we only need to determine the unknown coefficients from the right
hand side of the above formula. Since it is a linear polynomial in
three variables, we can easily know that the first two steps of the
$3D$ algorithm can be used to select four grid points (which are
not on the second layer down neighbor plane) which can be used to
uniquely determine the four unknown coefficients. Thus we have proved
that the algorithm can uniquely determine a quadratic polynomial in
3D. $\blacksquare$

Similarly as in $2D$, we can easily extend the algorithm to select
the stencil for higher order polynomials in 3D or extend it to higher
dimension.

\section{Application\label{sec:Examples}}

We now use the algorithm to solve an elliptic interface problem in
the Cartesian coordinate. For more detail, refer to \cite{WangSamulyakGuo2010,WangGlimmSamulyakJiaoDiao2013,SchBarCol2006}. 

\begin{figure}
\caption{\label{fig:Stencil-for-Partial-Cell}Stencil for the discretization
of a partial cell}

\begin{centering}
\includegraphics[scale=0.3]{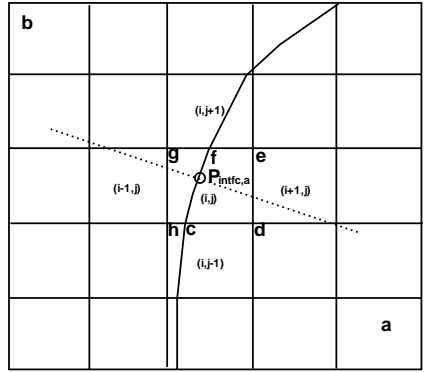}
\par\end{centering}

\end{figure}

An elliptic interface problem is a special elliptic problem with an
internal interface: \begin{equation}
\nabla\centerdot\frac{\nabla p}{\rho}=f\label{eq:elliptic}\end{equation}
where $\rho$ is a piecewise continuous function with jump across
the internal interface and $f$ is a given function which is continuous
inside each part of the domain. For the interior interface, we have
the following two jump conditions: \begin{equation}
[p]=J_{1}(\mathbf{x}),\label{eq:elliptic_potential_jump}\end{equation}
 \begin{equation}
[\frac{1}{\rho}\frac{\partial p}{\partial n}]=J_{2}(\mathbf{x}),\label{eq:elliptic_flux_jump}\end{equation}
 where $J_{1}$ and $J_{2}$ are given functions of the spatial variables
\cite{LeVequeLi1994}. The embedded boundary method \cite{WangSamulyakGuo2010,WangGlimmSamulyakJiaoDiao2013,SchBarCol2006}
uses a Cartesian mesh and the mesh cells are classified into four
types: external, internal, boundary, and partial cells depending on
whether there is an interior interface or exterior interface cutting
through the cells. The treatment of a boundary cell is similar to
that of an partial cells but much simpler. For an partial cell (see
Fig \ref{fig:Stencil-for-Partial-Cell}, cell $(i,j)$), it consists
of two parts separated by a cell interface. 

We store four unknowns at such a partial cell: two at the cell center
and two at the cell interface center for the two different components.
By using the elliptic equation \ref{eq:elliptic} for the two partial
parts of the cell and the two jump conditions, we could set up four
algebraic equations. To discretize equation \ref{eq:elliptic_flux_jump},
for example, we have \[
\left.\frac{1}{\rho}\frac{\partial p}{\partial n}\right|_{intfc,a}-\left.\frac{1}{\rho}\frac{\partial p}{\partial n}\right|_{intfc,b}=J_{2}(x)\]
where $\left.\frac{1}{\rho}\frac{\partial p}{\partial n}\right|_{intfc,a}$
is the flux calculated by polynomial interpolation selecting cells
belonging to component $a$ only. When we select cells for interpolation
using the $3D$ algorithm presented in the previous section, it is
best to try to find candidate cells in the normal direction.

Please refer to \cite{WangSamulyakGuo2010,WangGlimmSamulyakJiaoDiao2013,SchBarCol2006}
for more detail.

We use the method of manufactured solutions to verify our EBM implementation
for the elliptic interface problem in the Cartesian coordinate. The
computational domain is $[0,3]\times[0,3]\times[0,3]$. The interface
position is a sphere, given as\[
\sqrt{(x-1.5)^{2}+(y-1.5)^{2}+(z-1.5)^{2}}=1.\]
 We use $\rho_{d}=1$ inside the sphere and $\rho_{o}=0.001$ out
side. The exact solution is \[
p(x,y,z)=\begin{cases}
(x^{2}+y^{2}+z^{2})^{\frac{3}{2}}, & \sqrt{(x-1.5)^{2}+(y-1.5)^{2}+(z-1.5)^{2}}<1.\\
0.001\times(x^{2}+y^{2}+z^{2})^{\frac{3}{2}}+(1-\frac{1}{10})\times\frac{1}{8}, & other\end{cases}\]

We substitute the exact solution into the elliptic equation \ref{eq:elliptic}
to calculate the right hand side $f$. Table \ref{tab:Mesh-Convergence-Study-Elliptic}
shows the mesh convergence for the this problem. From the table, we
can see that the method is second order accurate in the $L_{\infty}$
norm. 

\begin{table}
\caption{\label{tab:Mesh-Convergence-Study-Elliptic}Mesh Convergence Study
for the Elliptic Interface Problem }

\centering{}\begin{tabular}{|c|c|}
\hline 
Mesh Size  & $L_{\infty}$ Error \tabularnewline
\hline
\hline 
10x10x10  & 0.19727211\tabularnewline
\hline 
20x20x20  & 0.05673376\tabularnewline
\hline 
40x40x40  & 0.01312201\tabularnewline
\hline 
80x80x80  & 0.00341078\tabularnewline
\hline
\end{tabular}
\end{table}

\section{Conclusion\label{sec:Conclusion}}

In this paper, we present a simple algorithm to select Cartesian grid
stencil for the multivariate interpolation in $2D$ and $3D$. The
correctness proof has been given using a constructive method. We show
its applicability by using this algorithm in the embedded boundary
method for solving the elliptic interface problem. 

It should be noted that the algorithm given is not unique. It is easy
to adopt the constructive method to create other algorithms for selecting
interpolation stencil for different needs. For example, we do not
need to restrict ourselves to select candidate grid points only in
the first and second layer neighbors. 

\bibliographystyle{plain}
\addcontentsline{toc}{section}{\refname}\bibliography{refs}

\end{document}